\documentclass[11pt,fleqn]{article}

\usepackage{rotating}
\usepackage{amsmath}   
\usepackage{amssymb}
\usepackage{latexsym}
\usepackage{graphicx}
\expandafter\ifx\csname
pdfoptionalwaysusepdfpagebox\endcsname\relax\else
\newtheorem{theorem}{Theorem}
\newtheorem{lemma}[theorem]{Lemma}

\newtheorem{obs}[theorem]{Observation}

\newtheorem{proposition}[theorem]{Proposition}
\newtheorem{ex}{Example}

\newcommand{\vx}{\boldsymbol x}
\newcommand{\vp}{\boldsymbol p}

\newcommand{\vu}{\boldsymbol u}
\newcommand{\vb}{\boldsymbol b}
\newcommand{\vd}{\boldsymbol d}
\newcommand{\zero}{\boldsymbol 0}

\newcommand{\vmu}{\boldsymbol \mu}
\newcommand{\vgamma}{\boldsymbol \gamma}
\newcommand{\va}{\boldsymbol a}
\newcommand{\mA}{\boldsymbol A}
\newcommand{\mB}{\boldsymbol B}
\newcommand{\mC}{\boldsymbol C}
\newcommand{\mI}{\boldsymbol I}
\newcommand{\mU}{\boldsymbol U}
\newcommand{\mD}{\boldsymbol D}
\newcommand{\mQ}{\boldsymbol Q}
\newcommand{\mM}{\boldsymbol M}
\newcommand{\mR}{\boldsymbol R}
\newcommand{\mS}{\boldsymbol S}
\newcommand{\mP}{\boldsymbol P}

\newcommand{\mT}{\boldsymbol T}
\newcommand{\vy}{\boldsymbol y}

\newcommand{\vt}{\boldsymbol t}
\newcommand{\vq}{\boldsymbol q}
\newcommand{\ve}{\boldsymbol e}

\newcommand{\vv}{\boldsymbol v}
\newcommand{\conv}{conv}

\newcommand{\qed}{\hfill$\Box$\vspace{\smallskipamount}}

\newcommand{\bbbz}{\mathbb{Z}}
\newcommand{\bbbr}{\mathbb{R}}

\title{Lattice based extended formulations for integer linear equality systems}

\author{Karen Aardal
\and Laurence A. Wolsey }

\date{February 12, 2007}

\begin{document}
\maketitle
\begin{abstract}
We study different extended formulations for the set $X =
\{\vx\in\mathbb{Z}^n \mid \mA\vx = \mA\vx^0\}$ in order to tackle
the feasibility problem for the set $X_+=X \cap \bbbz^n_+$. Here
the goal is not to find an improved polyhedral relaxation of
conv$(X_+)$, but rather to reformulate in such a way that the new
variables introduced provide good branching directions, and in
certain circumstances permit one to deduce rapidly that  the
instance is infeasible. For the case that $\mA$ has one row $\va$
we analyze the reformulations in more detail. In particular, we
determine the integer width of the extended formulations in the
direction of the last coordinate, and derive a lower bound on the
Frobenius number of $\va$. We also suggest how a decomposition of
the vector $\va$ can be obtained that will provide a useful
extended formulation. Our theoretical results are accompanied by a
small computational study.

\end{abstract}

\section{Introduction}
Recently, several hard integer programming feasibility problems
have been successfully tackled using a lattice reformulation
proposed by Aardal et al. \cite{AaHL98,AaHL00}. These problems
include various equality knapsacks of Cornu\'ejols et al.
\cite{CUWW}, Aardal and Lenstra \cite{AaL04}, market split
instances of Cornu\'ejols and Dawande \cite{CD99,AaBHLS}, and
financial planning instances of Louveaux and Wolsey \cite{LW00}.

The approach proposed earlier is to replace the feasibility
problem for the set $X=\{\vx \in \mathbb{Z}^n_+\mid
\mA\vx=\mA\vx^0\}$ by a feasibility problem over an extended
formulation for the set, namely the set $\{(\vx,\vmu)\in
\mathbb{Z}^n_+\times \mathbb{Z}^{n-m}\mid \vx=\vx^0+\mQ\vmu\}$
where $\mQ$ is a reduced basis of the lattice $\{\vy\in
\mathbb{Z}^n\mid \mA\vy=\zero\}$. It turns out that both a special
purpose branch-and-bound code and commercial mixed integer
programming solvers are then much more successful at solving the
test instances with the extended formulation than with the
original formulation. In related work Aardal and Lenstra
\cite{AaL04} have considered knapsack sets $\{\vx \in
\mathbb{Z}^n_+\mid \va\vx=b\}$ where $\va = M_1\vp^1+M_2\vp^2$,
$M_1,M_2\in\bbbz_{>0},\ \vp^1,\vp^2\in\bbbz^n$, where the vectors
$\vp^1$ and $\vp^2$ are short compared to the length of $\va$. For
such instances they have used the reduced basis approach to detect
a good direction for branching and to demonstrate infeasibility
for large values of $b$, and they theoretically derive  a strong
lower bound on the Frobenius number of $\va$ in the special case
when $M_2=1$ and $\vp^1\in\mathbb{Z}^n_{>0}$.

Here we pursue this lattice viewpoint.  In Section
\ref{sec_formulations} we show that the formulations presented
above are just two extremes of a family of extended formulations
of the original set $X$ taking the form $\{(\vx,\vmu)\in
\mathbb{Z}^n \times \mathbb{Z}^s\mid \mP\vx=\mP\vx^0+\mT\vmu\}$
with $s$ additional variables and $m+s$ constraints with $0 \leq
s\leq n-m$, where each row of $\mA$ can be expressed as an integer
multiple of the rows of $\mP$. Note that  whereas in polyhedral
combinatorics extended formulations are typically used to provide
better polyhedral approximations of $\conv(X)$, here the extended
formulations do not change the underlying polyhedron. However by
isolating good branching directions, the formulations are made
more effective for treatment by branch-and-bound.

In Section \ref{sec_K2}, we consider the special case of knapsack
sets when $m=1$ and $\va=M_1\vp^1+M_2\vp^2$. For the reformulation
with two equations and one additional variable ($s=1$), we
calculate the integer width of the underlying polyhedron in the
direction of the auxiliary variable $\mu$. This in turn leads us
to  a lower bound on the Frobenius number, simplifying and
generalizing the earlier proof from \cite{AaL04} cited above.

In Section \ref{sec_redbasis} we make precise how hidden structure
in $\mA$ can be detected by reduced basis calculations, and also
how this leads to interesting reformulations. In Section
\ref{sec_comp} we demonstrate computationally examples of hidden
structure, we test how effective the reformulations are in solving
feasibility problems, and we examine the quality of the lower
bounds on the Frobenius number.


\subsection{Preliminaries}\label{sec_prelim} In the following we
will frequently refer to lattices and reduced lattice bases. We
will define and describe these briefly here. For a more detailed
exposition, see for instance Cassels \cite{Cassels97}, Kannan
\cite{Kan87a}, Lenstra \cite{HWL00,HWL05}, and Lov\'asz
\cite{Lov86}.

Let $\vb^1,\dots,\vb^l$ be linearly independent vectors in
$\bbbr^n$. The set
\begin{displaymath}
L=\{\vx\in \bbbr^n \mid \vx=\sum_{j=1}^l\lambda_j\vb^j,\
\lambda_j\in\bbbz,\ 1\leq j\leq l\}
\end{displaymath}
is called a {\em lattice}. The set of vectors
$\{\vb^1,\dots,\vb^l\}$ is called a {\em lattice basis}. The
lattice generated by the basis $\mB$ is denoted by $L(\mB)$. If
$\mB$ is clear from the context we write just $L$.

Suppose $\mA$ is an integer $m\times n$ matrix of full row rank.
Then $N(\mA)=\{\vx\in\mathbb{R}^n\mid \mA\vx = \zero\}$, i.e.,
$N(\mA)$ is the null-space of $\mA$. We use the notation
$L(N(\mA))$ to denote the lattice that is generated by a basis of
all {\em integer vectors} in $N(\mA)$.

The {\em rank} of a given lattice $L$, rk $L$, is equal to the
dimension of the Euclidean vector space generated by a basis of
$L$.

Let $\mB=(\vb^1,\dots,\vb^l)$ be a basis for a given lattice $L$,
and let $\vb^{1*},\ldots,\vb^{n*}$ be the associated Gram-Schmidt
vectors, which are derived recursively as follows.
\begin{displaymath}
\vb^{1*}=\vb^1 \mbox{ and }
\vb^{j*}=\vb^j-\sum_{k=1}^{j-1}\mu_{jk}\vb^{k*},\quad 2\le j\le l,
\end{displaymath}
where
\begin{displaymath}
\mu_{jk}=((\vb^{j})^T\vb^{k*})/(||\vb^{k*}||^2),\quad 1\le k<j\le
l\,.
\end{displaymath}

The
basis $\mB=(\vb^1,\dots,\vb^l)$ of the lattice $L$ is
\textit{reduced in the sense of Lenstra-Lenstra-Lov\'asz}
\cite{LLL} if the following holds:
\begin{displaymath}
 |\mu_{jk}|\leq \frac{1}{2} \quad\mbox{ for }1\leq
k<j\leq l,
\end{displaymath}

\vspace*{-.4cm}
\begin{displaymath}
 ||\vb^{j*}+\mu_{j,j-1}\vb^{(j-1)*}||^2\geq
\frac{3}{4}||\vb^{(j-1)*}||^2 \quad\mbox{ for } 1<j\leq l.
\end{displaymath}
Such a basis can be found in polynomial time using the so-called
$LLL$ algorithm \cite{LLL}.

We will also need to find a basis $\mQ$ of $L(N(\mA))$ and a
vector $\vx^0$ satisfying $\mA\vx^0=\vb$, if such a vector exists.
They can be found by reducing the basis of the lattice defined by
the columns of the $(n+m+1)\times (n+1)$-matrix
\begin{displaymath}
\mB=\left(\begin{array}{cc}
\mI &  \zero\\
\zero & N_1\\
N_2\mA & -N_2\vb
\end{array}\right)\,,
\end{displaymath}
where $N_1$ and $N_2$ are appropriately chosen positive integers,
see \cite{AaHL00}. The first $n-m+1$ columns of the reduced basis
$\mB'$ will be
\begin{equation}\label{eq_AHLredbasis}
\left(\begin{array}{cc}
\mQ &  \vx^0\\
\zero & \pm N_1\\
\zero & \vb
\end{array}\right)\,.
\end{equation}
%
%
If there does not exists vector $\vx^0$ satisfying $\mA\vx^0=\vb$,
then element $b'_{n+1,n-m+1}$ of $\mB'$ will be equal to $\pm
kN_1$ for some integer $k>1$, indicating that there exists a
vector $\vx$ satisfying $\mA\vx=k\vb$.

The  \textit{Hermite Normal Form} of a matrix
$\mA\in\mathbb{Z}^{m\times n}$ of full row rank, HNF$(\mA)$, is
obtained by multiplying $\mA$ by an $n\times n$ unimodular matrix
$\mU$ to obtain the form $(\mD,\ \zero)$, where
$\mD\in\mathbb{Z}^{m\times m}$ is a nonsingular, nonnegative lower
triangular matrix with the unique row maximum along the diagonal.
The Hermite Normal Form can be calculated in polynomial time using
for instance the $LLL$ algorithm, see \cite{S86}, or that of
Kannan and Bachem \cite{KB79}.  
Note that
when $\mD \neq \mI$, the matrix $\mD^{-1}\mA$ is integral and
HNF$(\mD^{-1}\mA)=(\mI,\ \zero)$.

We will need to use one property concerning projections of
equality sets over the reals.

\begin{lemma}
\label{th_obs2} Let $\mC$ and $\mD$ be rational matrices of
dimension $m\times n$ and $m\times p$ respectively, and let $\vb$
be a rational $m$-vector. If
\begin{displaymath}
T=\{\vx\in\mathbb{R}^n \mid \mbox{ there exists
}\vy\in\mathbb{R}^p \mbox{ with } \mC\vx +\mD\vy = \vb\}\,,
\end{displaymath}
and if $\Delta$ is a basis of $N(\mD^T)$, then
\begin{displaymath}
T=\{\vx\in\mathbb{R}^n \mid \Delta^T\mC\vx = \Delta^T\vb\}\,.
\end{displaymath}
\end{lemma}

\begin{proof}
Multiplying by $\Delta^T$, it is immediate that
$\{\vx\in\mathbb{R}^n \mid \mbox{ there exists }\vy\in\mathbb{R}^p
\mbox{ with } \mC\vx +\mD\vy = \vb\}\ \subseteq
\{\vx\in\mathbb{R}^n \mid \Delta^T\mC\vx = \Delta^T\vb\}$ as
$\Delta^T\mD=0$. If the inclusion is strict, there exists
$\vx^*\in\mathbb{R}^n$ such that $\Delta^T\mC\vx^* = \Delta^T\vb$,
but $\{\vy\in \mathbb{R}^p \mid  \mD\vy =
\vb-\mC\vx^*\}=\emptyset$. This implies the existence of a vector
$\vu$ such that $\vu\mD=0$ and $\vu(\vb-\mC\vx^*) \neq 0$.
Howeever as $\Delta$ is a basis of $N(\mD^T)$, there exists $\vv$
such that $\vu=\vv\Delta^T$. But now
$\vu(\vb-\mC\vx^*)=\vv\Delta^T(\vb-\mC\vx^*)=0$, a contradiction.
\qed
\end{proof}

 The set of nonnegative (positive)
integers in $\mathbb{R}^n$ is denoted by $\mathbb{Z}^n_{+}$
($\mathbb{Z}^n_{>0}$). Similar notation is used for the
nonnegative real numbers.

\section{Formulations for equality constrained integer programs}\label{sec_formulations}
Our goal in this section is to derive reformulations of the
integer programming equality set
\begin{displaymath}
X=\{\vx\in\mathbb{Z}^n \mid \mA\vx=\vb\}
\end{displaymath}
and the associated affine set
\begin{displaymath}
Y=\{\vx\in\mathbb{R}^n \mid \mA\vx=\vb\},
\end{displaymath}
 where $\mA \in \mathbb{Z}^{m\times n}$ and $\vb\in\mathbb{Z}^n$.
 Later these will be used
 in testing feasibility or optimizing over the set $X^+=X \cap
\mathbb{Z}^n_+$.

 Throughout the paper we will assume that $X \neq \emptyset$. This can be tested by calculating HNF$(\mA)$.
  Note that for any point $\vx^0 \in X$, we can then write $X=\{\vx\in\mathbb{Z}^n \mid
  \mA\vx=\mA\vx^0\}$.

We will be dealing with the kernel lattice $L(N(\mA))=\{\vy \in
\mathbb{Z}^n\mid \mA\vy=\zero\}$. We assume that rank$(\mA)=m$, so
rk $(L(N(\mA)))=n-m$.

The following is immediate as $\vx \in X$ if and only if
$\vx-\vx^0 \in L(N(\mA))$.
\begin{obs}\label{th_obs1}
Let $\mQ \in \mathbb{Z}^{n\times(n-m)}$ be a lattice basis of
$L(N(\mA))$. Then
\begin{displaymath}
X=\{\vx^0\}+L(N(\mA))=\{\vx \mid \vx=\vx^0+\mQ\vmu, \ \vmu \in \mathbb{Z}^{n-m}\}.
\end{displaymath}
\end{obs}

Now we derive a larger family of reformulations of the set $X$.
Specifically suppose that two matrices $\mP \in
\mathbb{Z}^{(m+s)\times n}$ and $\mM \in \mathbb{Z}^{m \times
(m+s)}$ are given satisfying $\mA=\mM \mP $. Two natural questions
are:
\begin{enumerate}
\item Under what conditions does there exist an extended
formulation for $X$ of the form
$$X=\{x \in \mathbb{Z}^n \mid \mP\vx=\mP\vx^0+\mT\vmu, \vmu \in
\mathbb{Z}^s\}$$ where $\mT \in \mathbb{Z}^{(m+s)\times s}$?
\item
How should one choose the matrices $\mP,\mM,\mT$ so as to obtain
``interesting" extended formulations?
\end{enumerate}
We will use the following notation:
\begin{displaymath} U^{\mP,\mT}=\{(\vx,\vmu)\in \mathbb{Z}^n \times
\mathbb{Z}^s\mid \mP\vx=\mP\vx^0+\mT\vmu\} \ {\rm and \
}X^{\mP,\mT}={\rm proj}_{\vx}U^{\mP,\mT}.
\end{displaymath}

\vspace*{-15pt}
\begin{displaymath} V^{\mP,\mT}=\{(\vx,\vmu)\in \mathbb{R}^n \times
\mathbb{R}^s\mid \mP\vx=\mP\vx^0+\mT\vmu\},\ {\rm and}\
Y^{\mP,\mT}={\rm proj}_{\vx}V^{\mP,\mT}.
\end{displaymath}

Here, ${\rm proj}_{\vx}U^{\mP,\mT}=\{\vx\in\mathbb{Z}^n\mid \mbox{
there exists } \vmu\in\mathbb{Z}^s \mbox{ with }(\vx,\vmu)\in
U^{\mP,\mT}\}$. In this notation $X=X^{\mA,\zero}$ and Observation
\ref{th_obs1} states that $X=X^{\mI,\mQ}$ when $\mQ$ is a lattice
basis of $N(L(\mA))$.

In Propositions \ref{th_prop1} and \ref{th_prop2} below we address
the first question, and
 derive  conditions ensuring that
$X=X^{\mP,\mT}$.

\begin{proposition}\label{th_prop1}
Suppose that $\mA$, $\mP$ and $\mM$ are given with $\mP \in
\mathbb{Z}^{(m+s)\times n}$, $\mM \in \mathbb{Z}^{m\times (m+s)}$
and $\mA=\mM\mP$. Then
\begin{displaymath}X=X^{\mP,\mT}
\end{displaymath}
if $\mT\in \mathbb{Z}^{(m+s)\times s}$ is a lattice basis of
$L(N(\mM))$.
\end{proposition}

\begin{proof}
If $\vx \in X^{\mP,\mT}$, then there exists $\vmu$ such that
$(\vx,\vmu) \in U^{\mP,\mT}$. Left multiplying by $\mM$, we see
that $\mM\mP\vx=\mM\mP\vx^0+\mM\mT\vmu=\mM\mP\vx^0$, or
$\mA\vx=\mA\vx^0$. So $X^{\mP,\mT} \subseteq X$.

Conversely $\vx \in X$ implies that $\mA\vx=\mA\vx^0$, or that
$\mM\mP(\vx-\vx^0)=\zero$. As $\mP(\vx-\vx^0) \in
\mathbb{Z}^{m+s}$ and $\mT$ is a lattice basis of $L(N(\mM))$,
$\mP(\vx-\vx^0)=\mT\vmu$ for some $\vmu \in \mathbb{Z}^s$. Thus $X
\subseteq X^{\mP,\mT}$.\qed
\end{proof}

Note that a slightly stronger statement can be made.
\begin{proposition}\label{th_prop2}
 If the
Hermite Normal Form of $\mP$ is of the form $(\mI,\zero)$, then
\begin{displaymath}
X=X^{\mP,\mT}
\end{displaymath}
if and only if $\mT\in \mathbb{Z}^{(m+s)\times s}$ is a lattice
basis of $L(N(\mM))$.
\end{proposition}

\begin{proof}
Necessarily if $X=X^{\mP,\mT}$, then the columns of $\mT$ lie in
$L(N(\mM))$. If they do not form a lattice basis, there exists
$\vt^* \in L(N(\mM))$ that is not in the sublattice $L(\mT)$.
Also, as HNF$(\mP)=(\mI,\zero)$, there exists $\vx^*$ such that
$\mP\vx^*=\vt^*$. Now the point $\vy^*=\vx^0+\vx^*$ lies in $X$ as
$\mA\vy^*-\mA\vx^0=\mA\vx^*=\mM\mP\vx^*=\mM\vt^*=\zero$, but
$\mP(\vy^*-\vx^0)=\mP\vx^*=\vt^* \neq \mT\vmu$ for any $\vmu \in
\mathbb{Z}^s$, and thus $X\neq X^{\mP,\mT}$. \qed
\end{proof}

\vskip12pt Here we address the second question specifying one way
to choose appropriate matrices $\mP$ and $\mT$. The more
restrictive conditions imposed  will be used in Sections
\ref{sec_K2} and \ref{sec_redbasis} to make ``good" choices
 using reduced bases.

Consider a sublattice $L'$ of $L(N(\mA))$ of rank $0 \leq r \leq
n-m$ with lattice basis $\mR \in \mathbb{Z}^{n \times r}$, and
$\mQ=(\mR,\mS)$ an extension to a basis of $L(N(\mA))$. Let $\mP^T
\in \mathbb{Z}^{n \times n-r}$ be a lattice basis of
$L(N(\mR^T))=\{\vy \in \mathbb{Z}^n \mid \mR^T\vy=\zero\}$, and
$s=n-m-r$.

With $\mP$ constructed in this way, we have that:
\begin{enumerate}
\item[i)] $\mA=\mM\mP$ for some $m \times (m+s)$ integer matrix $\mM$.
$\mR \subseteq \mQ$ implies $L(N(\mQ^T)) \subseteq L(N(\mR^T))$,
and as $\mA^T \in L(N(\mQ^T))$ and $\mP^T \in L(N(\mR^T))$, the
claim follows.
\item[ii)] $\mT=\mP\mS$.
\end{enumerate}

\begin{theorem}\label{th_compform}
For any lattice basis $\mR$ of a sublattice of $L(N(\mA))$
\begin{eqnarray*}
&&X= X^{\mP,\mP\mS} \\
&&Y= Y^{\mP,\mP\mS}.
\end{eqnarray*}
\end{theorem}

\vspace{0.3cm}
\begin{proof}
  We show
that $X^{\mI,\mQ}\subseteq X^{\mP,\mP\mS} \subseteq X$. If $\vx\in
X^{\mI,\mQ}$, then there exists $\vmu$ such that $(\vx,\vmu) \in
U^{\mI,\mQ}$, and then multiplying through by $\mP$, we see that
$(\vx,\vmu)\in U^{\mP,\mP\mS}$ as $\mP\mR=\zero$ and thus $\vx\in
X^{\mP,\mP\mS}$. If $\vx\in X^{\mP,\mP\mS}$, then left multiplying
by $\mM$ gives $\mA\vx =\mA\vx^0=\vb$ as $\mM\mP\mS=\mA\mS=\zero$,
and thus $\vx\in X^{\mA,\zero}$. As $X^{\mI,\mQ}=X^{\mA,\zero}$,
equality holds for all $\mP$.

$Y^{\mI,\mQ}=Y^{\mP,\mP\mS}$ follows immediately from Lemma
\ref{th_obs2} as $\mP^T$ is a basis for $L(N(\mR^T))$, and
$Y^{\mI,\mQ}=Y^{\mA,\zero}$ follows from Lemma \ref{th_obs2} as
$\mA^T$ is a basis of $L(N(\mQ))$.\qed
\end{proof}

\begin{ex}\label{ex_ex1}{\rm
Consider a set $X=\{\vx \in \bbbz^5 | \mA\vx=\mA\vx^0\}$  of the
form

 \vskip10pt  $
\begin{array}{cccccccc}
x_1&-5x_2&-4x_3&+11x_4&+5x_5 &=& \va^1\vx^0\\
-13x_1&-2x_2&-12x_3&-11x_4&+x_5&=& \va^2\vx^0\\
&&x&\in&\mathbb{Z}^5.&&
\end{array}
 $

 \noindent Given
 \begin{displaymath}
 \mP= \left(
\begin{array}{rrrrr}
1&0&1&1&0\\
1&1&-1&1&0\\
0&-1&-1&2&1
\end{array}
\right) $ and  $\mM= \left(
\begin{array}{rrr}
1&0&5\\
-12&-1&1
\end{array}
\right),
\end{displaymath}
it is easily checked that $\mA=\mM\mP$.

\vskip10pt \noindent
Now as $\mT= \left(
\begin{array}{r}
-5\\
61\\
1
\end{array}
\right)$ is a basis of $L(N(\mM))$, it follows from Proposition
\ref{th_prop1} that
\begin{displaymath}
\begin{array}{cccccccc}
1x_1&&+1x_3&+1x_4&&=&\vp^1\vx^0&-5\mu\\
1x_1&+1x_2&-1x_3&+1x_4&&=&\vp^2\vx^0&+61\mu\\
&-1x_2&-1x_3&+2x_4&+1x_5&=&\vp^3\vx^0&+1\mu\\
&x&\in&\mathbb{Z}^5,&\mu&\in&\mathbb{Z}^1&
\end{array}
\end{displaymath}
is an extended formulation for $X$.

An alternative is to calculate a basis $\mQ$ of $L(N(\mA))$,
namely \begin{displaymath}
\mQ= \left(
\begin{array}{rrr}
1&-1&13\\
0&2&16\\
0&1&-25\\
-1&0&7\\
2&3&-22
\end{array}
\right)\,.
\end{displaymath}
Taking $\mR$ to consist of the first two columns, and $\mS$ to be
the last column, we have that $\mP^T$ from above is a basis of
$L(N(\mR^T))$ and $\mP\mS=\left(
\begin{array}{r}
-5\\
61\\
1
\end{array}
\right)$, so we arrive at the same extended formulation via the
formulation of Theorem \ref{th_compform}.
}\qed
\end{ex}

\section{Knapsack sets replaced by two equations}\label{sec_K2}

Here we analyze the case $m=1$ and $s=1$ in more detail.
Specifically we are interested in the sets $X$ and $Y$
 when
\begin{displaymath}
\va=M_1\vp^{1}+M_2\vp^{2}.
\end{displaymath}
We suppose that $a_j>0$ for all $1\le j \le n$, that
gcd$(a_1,\dots,a_n)=1$ and that $M_1,M_2>0$. Since
gcd$(a_1,\dots,a_n)=1$, it follows that gcd$(M_1,M_2)=1$, which
implies the existence of integers $\vq \in \mathbb{Z}^2$ with
$M_1q_1+M_2q_2=1$. In addition we add the condition that
HNF$(\mP)=(\mI,\zero)$. Note that if HNF$(\mP)=(\mD,\zero)$ with
$\mD \neq \mI$, then it suffices to multiply the equation system
by $\mD^{-1}$ and work with $\mP'=\mD^{-1}\mP$.

From Proposition \ref{th_prop1}, we know that $X^{\mP,\mT}=X$,
where $U^{\mP,\mT}$ is of the form
\begin{eqnarray*}
\vp^{1}\vx&=& \vp^{1}\vx^0 +M_2 \mu \\
\vp^{2}\vx&=& \vp^{2}\vx^0 -M_1 \mu \\
&&\vx \in \mathbb{Z}^n, \mu \in \mathbb{Z}^1.
\end{eqnarray*}
Note also that we can take $\vp^i\vx^0=q_ib$ for $i=1,2$ as
$M_1q_1+M_2q_2=1$. This follows because $\va\vx^0-\vb=
M_1(\vp^{1}\vx^0-q_1b)+M_2(\vp^{2}\vx^0-q_2b)=0$ and
gcd$(M_1,M_2)=1$.




\vskip10pt

In the rest of this section we add a nonnegativity constraint on
$\vx$. Thus we consider
\begin{eqnarray*}
X_+&=&X \cap \mathbb{Z}^n_+=\{\vx \in
\mathbb{Z}^n_+ \mid \va\vx=b\} \ {\rm and\ }\\
Y_+&=&Y \cap \mathbb{R}^n_+=\{\vx \in
\mathbb{R}^n_+ \mid \va\vx=b\}.
\end{eqnarray*}
We derive two results. The first concerns the width of the
polyhedron $V^{\mP,\mP\mS}_+$ in the direction $\mu$.
The second uses the width to derive a lower bound on the Frobenius
number, simplifying and generalizing the result of Aardal and
Lenstra \cite{AaL04} that is valid under the assumptions that
$\vp^1\in \mathbb{Z}_{>0}$ and $M_2=1$.

\subsection{The integer width}

The integer width of a rational polytope $P$ in the integer
direction $\vd$, $w_I(P,\vd)$, is defined as
\begin{displaymath}
w_I(P,\vd) = \left\lfloor\max\{\vd^T\vx \mid \vx\in
P\}\right\rfloor - \left\lceil\min\{\vd^T\vx \mid \vx\in
P\}\right\rceil +1\,,
\end{displaymath}
and is equal to the number of parallel lattice hyperplanes in
direction $\vd$ that are intersecting $P$.

Our goal is to calculate the integer width of
\begin{displaymath}
V^{\mP,\mP\mS}_+ =\{(\vx,\mu)\in\mathbb{R}^n_+\times
\mathbb{R}\mid
 \vp^1\vx-M_2\mu = q_1b,\ \vp^2\vx+M_1\mu = q_2b\}
\end{displaymath}
To do this,  let
\begin{equation}\label{eq_Vrhs1}
\bar{V}^{\mP,\mP\mS}_+ =\{(\vx,\mu)\in\mathbb{R}^n_+\times
\mathbb{R}\mid \vp^1\vx-M_2\mu = q_1,\ \vp^2\vx+M_1\mu = q_2\}
\end{equation} be the scaled down version of
this polyhedron with $b=1$.

Below we derive the values $\bar{z}=\max \{\mu \mid(\vx,\mu) \in
\bar{V}^{\mP,\mP\mS}_+\}$ and $\underline{z}=\min \{\mu
\mid(\vx,\mu) \in \bar{V}^{\mP,\mP\mS}_+\}$. Once we have the
values $\bar{z}$ and $\underline{z}$ it will be straightforward to
compute the width $w_I(V^{\mP,\mP\mS}_+,\ve^{n+1})$ as $\lfloor
b\bar{z}\rfloor - \lceil b\underline{z}\rceil +1$.


\begin{lemma} \label{lemma_width1} Consider formulation
$\bar{V}^{\mP,\mP\mS}_+$ (\ref{eq_Vrhs1}), and let $k={\rm
arg}\max\{i|p^1_i/a_i\}$ and $j={\rm arg}\min\{i|p^1_i/a_i\}$.
Then,
\begin{displaymath}
\bar{z}=\frac{p^1_k}{M_2a_k}-\frac{q_1}{M_2} \ {\rm and \ }
\underline{z}=\frac{p^1_j}{M_2a_j}-\frac{q_1}{M_2}.
\end{displaymath}
\end{lemma}

\begin{proof}
We consider the linear program
\begin{eqnarray}\nonumber
\bar{z}=\max \mu &&\\
\mbox{s.t. }\ \vp^{1}\vx-M_2\mu&=& q_1 \label{eq_pz}\\
\vp^{2}\vx+M_1\mu &=& q_2 \label{eq_rz}\\ \nonumber\vx \in
\mathbb{R}^n_+,&&\mu \in \mathbb{R}^1.
\end{eqnarray}

Let $\vgamma$ be the dual variables corresponding to constraints
(\ref{eq_pz})-(\ref{eq_rz}). The corresponding dual problem is:
\begin{eqnarray}\label{eq_obj}
\bar{z}=\min\ \ q_1\gamma_1+q_2\gamma_2 &&\\
\label{eq_dual1}
\mbox{s.t. }\ p^1_i\gamma_1+p_i^2\gamma_2 & \ge & 0,\quad 1\le i\le n,\\
\label{eq_dual2} -M_2\gamma_1 +M_1\gamma_2 &=& 1,\\ \nonumber
\vgamma \in \mathbb{R}^2 &&
\end{eqnarray}
From constraint (\ref{eq_dual2}) we obtain
\begin{equation}\label{eq_lambda}
\gamma_1 = \frac{M_1\gamma_2-1}{M_2}\,.
\end{equation}
Substituting for $\gamma_1$ in constraint (\ref{eq_dual1}) yields
\begin{displaymath}
p_i^1(\frac{M_1\gamma_2-1}{M_2})+p_i^2\gamma_2 \ge 0,\quad 1\le
i\le n\,.
\end{displaymath}
Rewriting gives $\gamma_2(p^2_i+p^1_i(M_1/M_2))-p^1_i/M_2\ge 0$
for $1\le i\le n$, which in turn yields
$\gamma_2\ge\frac{p^1_i}{a_i},\ 1\le i\le n$.  We now obtain
\begin{equation}\label{eq_mu}
\gamma_2=\frac{p^1_k}{a_k}\,.
\end{equation}
Finally, we substitute for $\vgamma$ in the dual objective
function (\ref{eq_obj}) using expressions (\ref{eq_lambda}) and
(\ref{eq_mu}) which yields the optimal dual objective value
\begin{displaymath}
\bar{z}=q_1\left(\frac{M_1(p^1_k/a_k)-1}{M_2}\right)+q_2
\frac{p^1_k}{a_k}=\frac{p^1_k}{a_k}\left(\frac{q_1 M_1+q_2
M_2}{M_2}\right)-\frac{q_1}{M_2}=\frac{p^1_k}{M_2a_k}-\frac{q_1}{M_2}\,.
\end{displaymath}

The calculation of $\underline{z}$ is almost identical.
 \qed
\end{proof}

\vspace{8pt}
Immediately we obtain the integer width.
\begin{theorem}\label{th_widthxz}
\begin{displaymath}
w_I(V_+^{\mP,\mP\mS},\ve^{n+1})
=\left\lfloor\frac{bp^1_k}{M_2a_k}-\frac{bq_1}{M_2}\right\rfloor -
\left\lceil\frac{bp^1_j}{M_2a_j}-\frac{bq_1}{M_2}\right\rceil
+1\,,
%
\end{displaymath}
where the indices $j$ and $k$ are defined as in Lemma
\ref{lemma_width1}.
\end{theorem}

Notice that the choice of a valid $\vq$ does not influence the
width $w(V^{\mP,\mP\mS}_+,\ve^{n+1})$. Suppose
$\vq'\in\mathbb{Z}^2$ satisfies $M_1q'_1+M_2q'_2=1$. The set of
all valid $(q_1,q_2)^T$ can be written as
\begin{displaymath}
\left(\begin{array}{c}
q_1\\
q_2\end{array}\right)= \left(\begin{array}{c}
q'_1\\
q'_2\end{array}\right)+\lambda \left(\begin{array}{r}
-M_2\\
M_1\end{array}\right),
\end{displaymath}
where $\lambda\in\mathbb{Z}$. Thus a different choice of $q_1$
just means that the whole interval $[\underline{z},\bar{z}]$ is
 shifted by an integer amount $\lambda$.

Observe also that given the choice of $\mR$, any lattice basis
$\mP$ of $L(N(\mR^T))$ yields the same width. The multipliers
$M_1,M_2$ as well as the values $q_1,q_2$ do however change with
different choices of $\mP$.

\begin{ex}\label{ex_ex2}{\rm
Consider the following input vector $\va$
\begin{displaymath}
\va =(12223, 12224, 36674, 61119, 85569).
\end{displaymath}
This is instance \texttt{cuww1} from Cornu\'ejols et al.
\cite{CUWW}.
With
\begin{displaymath}
\mP=\left(\begin{array}{rrrrr}
-1 & 0&2&-1&1\\
2 & 1&1&6&6
\end{array}\right)\,.
\end{displaymath}
and $(M_1,M_2)=(12225,12224)$, we have that $\va=\mM\mP$,
HNF$(\mP)=(\mI,\zero)$.

Now we apply Theorem \ref{th_widthxz}  with right-hand side
$b=89,643,481$, which is the Frobenius number of $\va$. We have
$j=1,\ k=3,\ q_1=1,\ q_2=-1$. We obtain
\begin{eqnarray*}
&&w_I(V^{\mP,\mP\mS}_+,\ve^{n+1})=\\
&& \left\lfloor
\frac{bp^1_k}{M_2a_k}-\frac{bq_1}{M_2}\right\rfloor - \left\lceil
\frac{bp^1_j}{M_2a_j}-\frac{bq_1}{M_2}\right\rceil +1 =
\left\lfloor\frac{bp^1_3}{M_2a_3}-\frac{bq_1}{M_2}\right\rfloor -
\left\lceil\frac{bp^1_1}{M_2a_1}-\frac{bq_1}{M_2}\right\rceil +1
=\\
&&\left\lfloor\frac{89643481}{12224}\left(\frac{2}{36674}-1\right)\right\rfloor
-
\left\lceil\frac{89643481}{12224}\left(\frac{-1}{12223}-1\right)\right\rceil
+1 = \\
&&\left\lfloor -7333.00003  \right\rfloor - \left\lceil
-7333.9999\right\rceil +1 = -7334+7333+1=0.
\end{eqnarray*}
 \vskip 10pt \noindent It follows that
$U_+^{\mP,\mP\mS}=\emptyset$. Applying branch-and-bound, and
branching first on the $\mu$ variable, this infeasibility would
immediately be apparent. This is not the case
 using branch and bound starting from the original formulation $X^{\mA,\zero}_+=\{x\in
\mathbb{Z}^n_+ \mid \va\vx=b\}$. In particular Cplex fails to
prove infeasibility within 500 million nodes. }
 \qed
\end{ex}

A natural question is whether the integer width differs if we use
a different member of the family of extended formulations.
Consider the sets
\begin{displaymath} V_+^{\mI,\mQ}=
\{(\vx,\vmu)\in \mathbb{R}_+^n \times \mathbb{R}^{s+r}\mid
\mI\vx=\mI\vx^0+\mS\vmu^S+\mR\vmu^R\}
\end{displaymath}
and
\begin{displaymath}
V_+^{\mP,\mP\mS}= \{(\vx,\vmu)\in \mathbb{R}^n_+ \times
\mathbb{R}^s\mid \mP\vx=\mP\vx^0+\mP\mS\vmu^S\}
\end{displaymath}
as described in Theorem \ref{th_compform}. Using Lemma
\ref{th_obs2}, we have that
\begin{displaymath}
{\rm proj}_{x,\mu^S}V_+^{\mI,\mQ}=V_+^{\mP,\mP\mS}.
\end{displaymath}
Thus we have
\begin{proposition}
$w_I(V_+^{\mP,\mP\mS},\vd)=w_I(V_+^{\mI,\mQ},(\vd,\zero))$, where
$\vd$ is any integer cost vector over the $(\vx,\vmu^S)$
variables.
\end{proposition}
In particular, when $m=1$, $s=1$, and $\vd$ is the unit vector
corresponding to the last column of $\mS$, then
$w_I(V_+^{\mP,\mP\mS},\ve^{n+1})=w_I(V_+^{\mI,\mQ},(\ve^{n+1},\zero))$.


\subsection{A lower bound on the Frobenius number}

The Frobenius number of $\va$, $F(\va)$, is the largest integer
value of $b$ such that $\va\vx=b$ does not have a nonnegative
integer solution. Wlog we choose $\vq$ such that $|q_1|\le M_2/2$.
So, if $|q'_1|>M_2/2$, we can determine new valid values of
$q_1,q_2$ such that $|q'_1|\le M_2/2$ by identifying an
appropriate value of $\lambda$. In this section we still assume
that HNF$(\mP)=(\mI,\zero)$.




\vspace{.3cm}
\begin{theorem}\label{th_froblb2}
Let $\va=M_1\vp^1+M_2\vp^2$ with $\va,M_1,M_2,\vp^1,\vp^2$
satisfying the assumptions given in the beginning of Section
\ref{sec_K2}. Moreover, let $\underline{z}$, $\bar{z}$ and the
indices $j$ and $k$ be as defined in Lemma \ref{lemma_width1}.

If $(-M_2/2)\le q_1 \le 0$ and
\begin{enumerate}
\item[1a)]  $\frac{p^1_j}{a_j}>q_1$ \item[2a)]
$\frac{p^1_k}{a_k}<M_2+q_1$ \item[3a)]
$(\frac{1-\bar{z}}{\bar{z}-\underline{z}})\underline{z}\not\in
\mathbb{Z}$
\end{enumerate}
then
\begin{equation}\label{eq_froblb2}
F(\va)\ge \frac{a_ja_k(M_2+q_1)-p^1_ka_j}{p^1_ka_j-p^1_ja_k}
-\frac{M_2}{\frac{p^1_j}{a_j}-q_1}\,.
\end{equation}
or if $0< q_1\le M_2/2$ and
\begin{enumerate}
\item[1b)]  $\frac{p^1_j}{a_j}>-M_2+q_1$
\item[2b)]  $\frac{p^1_k}{a_k}<q_1$
\item[3b)] $(\frac{1+\underline{z}}{\bar{z}-\underline{z}})\bar{z}\not\in
\mathbb{Z}$,
\end{enumerate}
then
\begin{displaymath}
F(\va)\ge \frac{a_ja_k(M_2-q_1)+p^1_ja_k}{p^1_ka_j-p^1_ja_k}
+\frac{M_2}{\frac{p^1_k}{a_k}-q_1}\,.
\end{displaymath}
\end{theorem}

\vspace{.3cm}
\begin{proof}

 We have already determined the width of $\bar{V}_+^{\mP,\mP\mS}$ in the direction
 of $\mu$ corresponding to $b=1$ in the proof of Lemma \ref{lemma_width1}.
 Specifically we have shown that $\mu$ lies in the interval
$[I_j,\ I_k]$, where
\begin{displaymath}
I_j:=\underline{z}=\frac{p^1_j}{M_2a_j}-\frac{q_1}{M_2} \mbox{ and
} I_k:=\bar{z}=\frac{p^1_k}{M_2a_k}-\frac{q_1}{M_2},
\end{displaymath} whose
width is
\begin{displaymath}
D:=I_k-I_j=\frac{a_jp^1_k-a_kp_j^1}{M_2a_ja_k}>0.
\end{displaymath}

Any integer right-hand side value $b=t$ for which the
corresponding interval $[tI_j,\ tI_k]$ does not contain an integer
is a lower bound on the Frobenius number $F(\va)$. Below we will
show that
\begin{displaymath}
t\ge \frac{a_ja_k(M_2+q_1)-p^1_ka_j}{p^1_ka_j-p^1_ja_k}
-\frac{M_2}{\frac{p^1_j}{a_j}-q_1}
\end{displaymath}
is such a lower bound in the case that $-M_2/2\le q_1<0$. A sketch
of the proof for the case $0< q_1\le M_2/2$ is given in Appendix
1.

If $q_1\le 0$,  Assumptions 1a and 2a imply that $0<I_j<I_k<1$.
Moreover, since $q_1\ge -M_2/2$ we obtain $I_k\le
p^1_k/(M_2a_2)+1/2$. Let $s:=\frac{1-I_k}{D}$. Notice that
$1-I_k>0$ since $I_k<1$. The interval $[sI_j,\ sI_k]$ has length
$1-I_k$. Notice that $sI_j\not\in\mathbb{Z}$ due to Assumption 3
of the theorem. Define $\ell :=\lfloor sI_j\rfloor$ and
$s':=\ell/I_j$. The number $s'$ satisfies $s-\frac{1}{I_j}< s'<
s$, and yields the interval $[I'_j,I'_k]:=[s'I_j,\ s'I_k]$, with
$I'_j$ integral. The length of $[I'_j,I'_k]$ is less than the
length $1-I_k$ of $[sI_j,\ sI_k]$. Therefore, $[I'_j,I'_k+I_k]$
has length less than 1, and since $I'_j$ is integral it follows
that $(I'_j,I'_k+I_k]$ does not contain an integer.

Now, define $s^{\ast}:=\lfloor s'\rfloor+1$ and the interval $[I^{\ast}_j,I^{\ast}_k]:=
[s^{\ast}I_j,\ s^{\ast}I_k]$. We have $I'_j<I^{\ast}_j\le I'_j+I_j$
and $I'_k<I^{\ast}_k\le I'_k+I_k$. The result that $[I^{\ast}_j,\ I^{\ast}_k]$ does
not contain an integer follows from the observation that $(I'_j,I'_k+I_k]$
does not contain an integer.

We finally observe that
\begin{displaymath}
s^{\ast}=\lfloor s'\rfloor +1 > \lfloor s-\frac{1}{I_j}\rfloor +1 \ge
s-\frac{1}{I_j} -1 +1 = s-\frac{1}{I_j}\,,
\end{displaymath}
so we can conclude that
$s-\frac{1}{I_j}=\frac{1-I_k}{D}-\frac{1}{I_j}$ yields a lower
bound on the Frobenius number $F(\va)$. Rewriting
$\frac{1-I_k}{D}-\frac{1}{I_j}$ results in the expression
\begin{displaymath}
\frac{1-I_k}{D}-\frac{1}{I_j} = \frac{1-\frac{p^1_k}{M_2a_k}+\frac{q_1}{M_2}}
{\frac{a_jp^1_k-a_kp^1_j}{M_2a_ja_k}}-\frac{1}{\frac{p^1_j}{M_2a_j}-\frac{q_1}{M_2}}=
\frac{a_ja_k(M_2+q_1)-p^1_ka_j}{p^1_ka_j-p^1_ja_k}
-\frac{M_2}{\frac{p^1_j}{a_j}-q_1}\,.
\end{displaymath}

\qed
\end{proof}


We notice the similarity with the expression for the lower bound
on the Frobenius number derived by Aardal and Lenstra \cite{AaL04}
for the case that $M_2=1$ and $\vp^1\in\mathbb{Z}^n_{>0}$. If we
set $M_2=1$ and $q=0$ in Expression (\ref{eq_froblb2}) we obtain
\begin{displaymath}
\frac{a_ja_k-p^1_ka_j}{p^1_ka_j-p^1_ja_k}-\frac{a_j}{p^1_j}.
\end{displaymath}
The only difference in the two expressions is in the numerator of the first
term, where we have $p^1_ka_j$ instead of $2p^1_j a_k$ in \cite{AaL04}.
This is a result of a different choice of the number $s$ in the proof. In \cite{AaL04}
$s$ was chosen as $s=(1-2I_j)/D$ under a constraint on the relationship
between $I_j$ and $I_k$.

\section{Computation}
\subsection{Using reduced bases to find structure}\label{sec_redbasis}
Aardal and
Lenstra \cite{AaL04} considered knapsack instances in which the
input vector $\va$ can be decomposed as $\va=M_1\vp^1+M_2\vp^2$,
and used the reformulation earlier suggested in \cite{AaHL00}, cf.
formulation $X^{I,Q}$:
\begin{displaymath}
\vx=\vx^0+\mQ\vmu\,,
\end{displaymath}
where $\vx^0$ and $\mQ$ is as described in Section
\ref{sec_formulations}. They observed that if $\va$ is long,
$\vp^1$ and $\vp^2$ are short compared to $\va$, and if the basis
$\mQ$ is reduced, then the first $n-2$ basis vectors of $\mQ$ are
short, and the last, $(n-1)$st basis vector is long. This can be
explained as follows. First we observe that the orthogonal
complement of the plane spanned by the vectors $\vp^1$ and $\vp^2$
is equal to $N(\mP)\subset N(\va)$. If $\vp^1$ and $\vp^2$ are
short, the lattice $L(N(\mP))$ contains short vectors yielding a
relatively small lattice determinant. The rank of $L(N(\mP))$,
which is a sublattice of $L(N(\va))$, is just one less than the
rank of $L(N(\va))$. Moreover, the determinant of the lattice
$L(N(\va))$ is equal to the length of the vector $\va$. So, the
large value of $d(L(N(\va)))$ mainly has to be contributed by the
basis vector that is in $L(N(\va))$ but not in $L(N(\mP))$. Since
basis reduction orders the basis vectors in nondecreasing order of
length, up to a multiplier, this basis vector is the last one in a
reduced basis.

In general, we can derive a suitable decomposition of $\mA$, in
case no such decomposition is known a priori, by using the
following algorithm.

\begin{enumerate}
\item[i)] Derive a reduced  basis $\mQ$ of $L(N(\mA))$, see
(\ref{eq_AHLredbasis}) in Section \ref{sec_prelim}.
\item[ii)]
Suppose $\mQ$ consists of $s$ long vectors and $r=n-m-s$ short
ones. How to define ``long'' and ``short'' is up to the user. (If
all vectors of $\mQ$ are of approximately the same length we set
$s=n-m$).  We define $\mR$ to be the set of short vectors of $\mQ$
and $\mS$ to be the set of long ones.
\item[iii)] Find a reduced
basis $\mP^T$ of $L(N(\mR^T))$.
\item[iv)] Solve the system of equations $\mM\mP=\mA$, $\mM \in
\mathbb{Z}^{m \times (m+s)}$ to find the matrix of multipliers
$\mM$.
\end{enumerate}

\begin{ex}\label{ex_ex3}{\rm
We consider the same instance as in Example \ref{ex_ex2} with
\begin{displaymath}
\va =(12223, 12224, 36674, 61119, 85569).
\end{displaymath}
Here we show how its hidden structure can be uncovered.

 A reduced basis of $L(N(\va))$ is equal to
\begin{displaymath}
\mQ = \left(\begin{array}{rrrr}
0 & -3 & -1 & 2059\\
1 & 1 & -3 & 157\\
-1 & -1 & -1 & -3336\\
-1 & 1 & 0 & 2687\\
1 & 0 & 1 & -806
\end{array}\right)\,.
\end{displaymath}
Here we observe that the last column of the reduced basis $\mQ$ is
much longer than the other columns. Taking $r=3$ and $s=1$,  $\mR$
will consist of the first three columns of $\mQ$, and $\mS$ will
consist of the last column of $\mQ$. A reduced basis $\mP^T$ for
the lattice $L(N(\mR^T))$ is
\begin{displaymath}
\mP^T=\left(\begin{array}{rr}
-1 & 2\\
0 & 1\\
2 & 1\\
-1 & 6\\
1 & 6
\end{array}\right)\,.
\end{displaymath}
The vector $(M_1,M_2)=(12225,12224)$ solves $\mM\mP=\va$. We can
now write $\va=12225\vp^1 +12224\vp^2$, with $\vp^1$ being the
first row of $\mP$, and $\vp^2$ being the second row of $\mP$.
Note that the matrix $\mP$ obtained is not unique. In Example 2 a
closely related but different basis of  $L(N(\mR^T))$ is
considered.}\qed
\end{ex}

\subsection{Feasibility testing and quality of the Frobenius bound}\label{sec_comp} \noindent

We tested the quality of the extended formulations for different
choices of $s$ on some instances of integer equality knapsacks and
the Cornu\'ejols-Dawande market split problem. For all instances
of both problem types we compute a reduced basis $\mQ$ and a
vector $\vx^0$ as described in Section \ref{sec_prelim}, and
 derive matrices $\mP$ and $\mM$ as described in Section
\ref{sec_redbasis}.

The integer knapsack instances were taken from Aardal and Lenstra
\cite{AaL04}. Instances prob1--4 are  such that the vector $\va$
decomposes with short $\vp^1,\ \vp^2$, whereas for the instances
prob11--14 the $\va$-coefficients, randomly generated from
$U[10000,150000]$, are of the same size on average as in prob1--4.
Instances prob11--14 have no apparent structure, and the columns
of a reduced basis $\mQ$ of $L(N(\mA))$ are of approximately the
same length. We use the Frobenius number of the vector $\va$ as
right-hand side coefficient for all knapsack instances. Instances
prob1--4 have 8 variables and prob10--14 have 10 variables. For
details of the instances, see \cite{AaL04}.

The market split instances \cite{CD99} are multiple row equality
knapsack problems in $\{0,1\}$-variables with $m$ rows and
$n=10(m-1)$ variables. The elements of $\va^i$ for each row $i$
are generated randomly from $U[0,99]$, and the right-hand side
coefficients are calculated as $b_i=\lfloor (\sum_{j=1}^n
a^i_{j})/2\rfloor$. We generated two sets of market split
instances with 4 constraints and 30 variables, and 5 constraints
and 40 variables respectively.

\begin{table}[h]
\caption{The number of branch-and-bound nodes needed to solve the
various reformulations for the knapsack
instances.}\label{table_nonodesknpsk}
\renewcommand{\arraystretch}{1.0}
\setlength\tabcolsep{5pt}
\begin{center}
\begin{tabular}{|l||r|r|r|r|r|r|r|r|r|r|}\hline
           &                            &    &
           \multicolumn{8}{c|}{$U_+^{\mP,\mP\mS}$}\\
Instance   & \multicolumn{1}{c|}{orig} & AHL &  $s=1$ & $s=2$ &
$s=3$ & $s=4$ & $s=5$ & $s=6$ & $s=7$& $s=9$\\ \hline
prob1 & $>100$ mill. &  1 & 59 & 15 & 3 & 3 & 1 & 1 & 1 &--\\
prob2 & $>100$ mill. &  3 & 23  & 7  & 3 & 1 & 1& 1 & 1 &--\\
prob3 & $>100$ mill. & 13 & 37 & 29  & 5 & 7 & 11 &9 & 5 &--\\
prob4 & $>100$ mill. &  3 & 13 & 5 & 1 & 1 & 1 & 1 & 1 &--\\
\hline
prob11 & 100,943 & 61 &  2237 & 7683 & 317 & 89 & 51 &69 & 49&  61\\
prob12 & 160,783 & 93 &  10,981 & 1105 & 967 & 523 & 179& 105 &117 &71\\
prob13 & 188,595 & 91 &  10,205 & 12,261 & 239 & 321 & 35& 57& 39 &59 \\
prob14 & 140,301 & 87 & 2443 & 627 & 689 & 389 & 283 & 115 &105 & 87\\
\hline
\end{tabular}
\end{center}
\end{table}

In Tables \ref{table_nonodesknpsk}--\ref{table_nonodesCD5} we
report on the number of nodes used by the integer programming
solver Xpress Version 16.01.01 \cite{Xpress}  to solve the various
reformulations. Column ``orig'' refers to the original formulation
in $\vx$-variables. Column ``AHL'' refers to the
Aardal-Hurkens-Lenstra lattice reformulation in which the
$\vx$-variables have been removed from the formulation, i.e., the
formulation $\{\vmu\in\mathbb{Z}^{n-1}\mid \mQ\vmu\ge -\vx^0\}$ in
the knapsack case and the formulation
$\{\vmu\in\mathbb{Z}^{n-m}\mid -\vx^0\le \mQ\vmu\le 1-\vx^0\}$ in
the market split case. For formulations $U_+^{\mP,\mP\mS}$ we
report on results for different values of $s$. Notice that the
formulations AHL and $U_+^{\mP,\mP\mS}$ for $s=n-m$ are
mathematically equivalent, but the $U_+^{\mP,\mP\mS}$-formulations
contain the $\vx$-variables with the identity matrix as
coefficients. Since the solver reacts differently to the presence
of the redundant $\vx$-variables, this leads to slight deviations
in the number of enumeration nodes needed.

Instances prob1--4, which decompose in short $\vp^1,\vp^2$, are
very difficult to tackle with branch-and-bound applied to the
original formulation. The Frobenius numbers for these instances
are also large, see Table \ref{table_froblb}. None of the
instances could be solved within 100 million nodes. As could be
expected, the $U_+^{\mP,\mP\mS}$-formulation with $s=1$, which is
a formulation with the $\vx$-variables and one variable $\mu$, is
easy to solve and comparable to the AHL-formulation. In contrast,
instances prob11--14 are solvable using the original formulation,
mainly due to the smaller value of the right-hand side
coefficients. Here, one could expect that we would need to set
$s=n-m$ to see a noticeable improvement compared to the original
formulation, but in fact even taking $s=1$ reduces the number of
enumeration nodes by at least an order of magnitude, and with $s$
around 5 we obtain results comparable to those obtained with the
AHL-formulation.

\begin{table}[h]
\caption{The number of branch-and-bound nodes needed to solve the
various reformulations: CD-instances $4\times
30$.}\label{table_nonodesCD4}
\renewcommand{\arraystretch}{1.0}
\setlength\tabcolsep{5pt}
\begin{center}
\begin{tabular}{|l||r|r||r|r|r|r|r|r|}\hline
           &                            &    &
           \multicolumn{6}{c|}{$U_+^{\mP,\mP\mS}$-ref}\\
Instance   & \multicolumn{1}{c|}{orig} & AHL &  $s=1$ &
$s=5$&$s=10$ & $s=15$ & $s=20$ & $s=26$\\ \hline
$4\times 30$\_1 & 157,569  & 281 & 124,695 & 71641 & 8033 & 1397 & 1021 & 607\\
$4\times 30$\_2 &  169,455 & 167 & 154,505&51989 & 3794 & 1487 & 610 & 535  \\
$4\times 30$\_3 &  209,741 & 325 & 178,697&181,373 & 32,367 & 1831 & 1025 & 845  \\
$4\times 30$\_4 &  202,513 & 199  &156,047 & 4685& 3583 & 829 & 493&9527 \\
$4\times 30$\_5 &  115,173 & 311 & 73,151&17,201 & 1197 & 391 & 353 & 3135  \\
\hline
\end{tabular}
\end{center}
\end{table}

\begin{table}[h]
\caption{The number of branch-and-bound nodes needed to solve the
various reformulations: CD-instances $5\times
40$.}\label{table_nonodesCD5}
\renewcommand{\arraystretch}{1.0}
\setlength\tabcolsep{5pt}
\begin{center}
\begin{tabular}{|l||r|r||r|r|r|r|r|}\hline
           &                            &    &
           \multicolumn{5}{c|}{$U_+^{\mP,\mP\mS}$-ref}\\
Instance   & \multicolumn{1}{c|}{orig} & AHL &
\multicolumn{1}{c|}{$s=5$} & \multicolumn{1}{c|}{$s=10$} &
\multicolumn{1}{c|}{$s=20$} & \multicolumn{1}{c|}{$s=30$} & \multicolumn{1}{c|}{$s=35$} \\
\hline
$5\times 40$\_1 & $>10,000,000$  & 5873  & $>10,000,000$ & 3,144,737 & 160,701 & 32,507 & 32,099   \\
$5\times 40$\_2 & $>10,000,000$  &  1643 & $>10,000,000$ & 2,821,042 & 128,707 & 30,302 & 12,734  \\
$5\times 40$\_3 & $>10,000,000$  & 7349  & $>10,000,000$ & 8,264,955 & 86,483 & 28,491 & 25,541  \\
$5\times 40$\_4 & $>10,000,000$  & 6870  & $>10,000,000$ &1,854,280  & 70,949 & 19,616 & 16,557 \\
$5\times 40$\_5 & $>10,000,000$  &  6651 & $>10,000,000$ & 7,805,023 &1,107,713  & 35,989 & 36,897 \\
\hline
\end{tabular}
\end{center}
\end{table}

For the market split instances, which have no clear structure of
the $\mQ$-matrix, we notice similar results to those obtained for
the knapsack instances prob11--14. The algorithm of Section
\ref{sec_redbasis} prescribes $s=n-m$ for these types of
instances. The computational results suggest that smaller values
of $s$ already yield significant computational improvement.

In Table \ref{table_froblb} we report on the value of the
Frobenius number as well as the value produced by the lower bound
given in Theorem \ref{th_froblb2}. For instances prob1--4 the
lower bound is of the same order of magnitude as the Frobenius
number, whereas for instances prob11--14 the bound is off by an
order of magnitude. The bound might be improved by a different
choice of the value $s$ in the proof of the theorem.

\begin{table}[h]
\caption{The value of the lower bound of the Frobenius
number.}\label{table_froblb}
\renewcommand{\arraystretch}{1.0}
\setlength\tabcolsep{5pt}
\begin{center}
\begin{tabular}{|l||r|r|}\hline
Instance   & \multicolumn{1}{c|}{$F(\va)$} & lower bound on
$F(\va)$\\
\hline
prob1 & 33,367,335 &  26,061,675 \\
prob2 & 14,215,206 & 10,894,273   \\
prob3 & 58,424,799 & 31,510,625  \\
prob4 & 60,575,665 & 56,668,034   \\
\hline
prob11 & 577,134 & 98,774  \\
prob12 & 944,183 & 113,114   \\
prob13 & 765,260 & 67,752    \\
prob14 & 680,230 & 60,476  \\
\hline
\end{tabular}
\end{center}
\end{table}


\newpage
\section*{Acknowledgements} This work was
partly carried out within the framework of ADONET, a European
network in Algorithmic Discrete Optimization, contract no.
MRTN-CT-2003-504438.
The first author is  financed in part by the Dutch BSIK/BRICKS
project. The research was carried out in part while the second
author visited CWI, Amsterdam with the support of the NWO visitor
grant number B 61-556.

\section*{Appendix 1}

\noindent \textbf{Proof of Theorem 6 for the case $0<q_1\le
(M_2/2)$.}\\ If $0< q_1\le (M_2/2)$,  Assumptions 1b and 2b imply
that $-1<I_j<I_k<0$, so the interval $[I_j,\ I_k]$ does not
contain an integer. In addition, $I_j\ge p^1_j/(M_2a_j)-1/2$.

Let $s:=\frac{1+I_j}{D}$. The length of the interval $[sI_j,\
sI_k]$ is equal to $1+I_j$, and since $-1<I_j<0$ we have that
$0<1+I_j<1$.

Notice that $sI_k\not\in\mathbb{Z}$ due to Assumption 3b of the
theorem. Define $\ell :=\lceil sI_k\rceil$ and $s':=\ell/I_k$. The
number $s'$ satisfies $s+\frac{1}{I_k}< s'< s$, and yields the
interval $[I'_j,I'_k]:=[s'I_j,\ s'I_k]$, with $I'_k$ integral. The
length of $[I'_j,I'_k]$ is less than the length $1+I_j$ of
$[sI_j,\ sI_k]$. Therefore, $[I'_j+I_j,I'_k]$ has length less than
1, and since $I'_k$ is integral it follows that $[I'_j+I_j,I'_k)$
does not contain an integer.

Now, define $s^{\ast}:=\lfloor s'\rfloor+1$ and the interval
$[I^{\ast}_j,I^{\ast}_k]:= [s^{\ast}I_j,\ s^{\ast}I_k]$. We have
$I'_j+I_j\le I^{\ast}_j< I'_j$ and $I'_k+I_k\le I^{\ast}_k< I'_k$.
The result that $[I^{\ast}_j,\ I^{\ast}_k]$ does not contain an
integer follows from the observation that $[I'_j+I_j,I'_k)$ does
not contain an integer.

We finally observe that
\begin{displaymath}
s^{\ast}=\lfloor s'\rfloor +1 > \lfloor s+\frac{1}{I_k}\rfloor +1
\ge s+\frac{1}{I_k} -1 +1 = s+\frac{1}{I_k}\,,
\end{displaymath}
so we can conclude that
$s+\frac{1}{I_k}=\frac{1+I_j}{D}+\frac{1}{I_k}$ yields a lower
bound on the Frobenius number $F(\va)$. Rewriting
$\frac{1+I_j}{D}+\frac{1}{I_k}$ results in the expression
\begin{displaymath}
\frac{1+I_k}{D}+\frac{1}{I_k} =
\frac{1+\frac{p^1_j}{M_2a_j}-\frac{q_1}{M_2}}
{\frac{a_jp^1_k-a_kp^1_j}{M_2a_ja_k}}+\frac{1}{\frac{p^1_k}{M_2a_k}-\frac{q_1}{M_2}}=
\frac{a_ja_k(M_2-q_1)+p^1_ja_k}{p^1_ka_j-p^1_ja_k}
+\frac{M_2}{\frac{p^1_k}{a_k}-q_1}\,.
\end{displaymath}

\end{document}